\documentclass[11pt,a4paper]{article}
\usepackage[T1]{fontenc}
\usepackage[utf8]{inputenc}
\usepackage{lmodern,amsmath,amssymb,amsthm,mathtools}
\usepackage[margin=28mm]{geometry}
\usepackage[expansion=false]{microtype}
\usepackage{indentfirst}
\usepackage{xcolor}
\definecolor{linkblue}{RGB}{25,61,127}
\usepackage[unicode,bookmarksopen=true,bookmarksnumbered=true,
  colorlinks=true,linkcolor=linkblue,citecolor=linkblue,urlcolor=linkblue,
  linktoc=all,pdfborder={0 0 0}]{hyperref}
\usepackage{bookmark}
\hypersetup{
  pdftitle={Bernstein--Sato polynomials of generic determinantal ideals},
  pdfauthor={Yifan Chen and Huaiqing Zuo},
  pdfsubject={Roots, multiplicities, and the strong motivic monodromy conjecture},
}
\numberwithin{equation}{section}
\newtheorem{maintheorem}{Theorem}

\newtheorem{theorem}{Theorem}[section]
\newtheorem{proposition}[theorem]{Proposition}
\newtheorem{lemma}[theorem]{Lemma}
\newtheorem{corollary}[theorem]{Corollary}
\theoremstyle{definition}

\newtheorem{definition}[theorem]{Definition}
\newtheorem{remark}[theorem]{Remark}
\newcommand{\C}{\mathbb C}
\newcommand{\Z}{\mathbb Z}
\newcommand{\N}{\mathbb N}
\newcommand{\g}{\mathfrak g}
\DeclareMathOperator{\Ann}{Ann}
\DeclareMathOperator{\Sym}{Sym}
\DeclareMathOperator{\Span}{span}
\DeclareMathOperator{\gr}{gr}
\DeclareMathOperator{\mult}{mult}
\DeclareMathOperator{\Mat}{Mat}

\title{Capelli Operators and the Bernstein-Sato polynomials of Minors}
\author{
  Yifan Chen\thanks{
    Department of Mathematical Sciences,
    Tsinghua University,
    Beijing 100084, P. R. China.
    Email: \texttt{c-yf20@tsinghua.org.cn}.
  }
  and
  Huaiqing Zuo\thanks{
    Department of Mathematical Sciences,
    Tsinghua University,
    Beijing 100084, P. R. China.
    Email: \texttt{hqzuo@mail.tsinghua.edu.cn}.
  }
}

\date{}
\begin{document}
\maketitle
\begin{abstract}
For every ideal of minors of a generic complex matrix, we determine all roots and their multiplicities in the Bernstein-Sato polynomial. These roots lie in explicit rational intervals indexed by the rank strata. Every nonintegral root is simple, whereas the multiplicity of each integral root equals the number of labels producing that root. We deduce the strong motivic monodromy conjecture, including the inequality relating pole orders to root multiplicities; for minors of size $r$, the largest possible root multiplicity over all matrix dimensions is $\lfloor\sqrt{4r-1}\rfloor$.
	

\end{abstract}

\begingroup
\small
\tableofcontents
\endgroup

\section{Introduction}\label{sec:introduction}

For a nonzero polynomial $f$ over $\C$, its Bernstein--Sato polynomial
is the monic polynomial $b_f(s)$ of least degree for which
\[
P(s)f^{s+1}=b_f(s)f^s
\]
with $P(s)$ a polynomial in $s$ whose coefficients are polynomial
differential operators. This equation arose in the work of Bernstein
and Sato on differential equations and complex powers. Bernstein's
existence theorem gives meromorphic continuation of distributions
depending on a complex parameter \cite{Bernstein}. Kashiwara proved
the rationality of the roots \cite{Kashiwara}. The resulting invariant
connects differential operators, singularities, and local monodromy.

Budur, Musta\c t\u a, and Saito extended the construction to arbitrary
ideals \cite{BMS}. For generators $f_1,\ldots,f_a$, their definition
uses the graph module and the Euler operator in its $a$ auxiliary
directions. The resulting single-variable polynomial is independent
of the chosen generators. Its largest root is the negative of the
log canonical threshold; more generally, the jumping coefficients in
$[\operatorname{lct}(I),\operatorname{lct}(I)+1)$ give roots after changing sign
\cite[Theorem~2]{BMS}. These
connections make ideals of minors a natural test family: their
rank strata, resolutions, and symmetry actions are explicit, while
their equations usually have many generators.

Let $X=(x_{ab})$ be an $m\times n$ matrix of independent variables over
$\C$, where $1\le r\le m\le n$. Write $I_r=I_r^{m,n}$ for the ordinary
ideal generated by all $r\times r$ minors in $S=\C[x_{ab}]$. We use the monic,
unshifted, single-variable ideal polynomial $b_{I_r}(s)$ of \cite{BMS};
its graph-module definition is recalled in Section~\ref{subsec:graph}.
The variety $V(I_r)$ consists of matrices of rank less than $r$.
Near a matrix of rank $i<r$, its transverse equations are the
$(r-i)$-minors of a generic $(m-i)\times(n-i)$ matrix. This recursive
geometry governs several of its singularity invariants.

Docampo studied the orbits in the arc spaces of generic determinantal
varieties and used them to compute jet-scheme components, log
canonical thresholds, and topological zeta functions \cite{Docampo}.
In the present notation, the threshold is
\[
\operatorname{lct}(I_r)
=\min_{0\le i<r}\frac{(m-i)(n-i)}{r-i}.
\]
In this integrability description, $f_j$ ranges over the $r$-minors, and the log canonical threshold
is the supremum of real $c$ for
which $(\sum_j|f_j|^2)^{-c}$ is locally integrable at every point.
Docampo also records Johnson's earlier computation of multiplier
ideals using resolutions by complete collineations. Thus the
largest Bernstein--Sato root was accessible geometrically before a
formula for the full polynomial.

The full calculation for generic minors was posed explicitly in
Budur's lecture notes \cite[Section~3.11]{BudurNotes}.
L\H orincz, Raicu, Walther, and Weyman solved the maximal-minor case
\cite{LRWW}, obtaining
\[
b_{I_m^{m,n}}(s)=\prod_{j=n-m+1}^{n}(s+j).
\]
Their method combines invariant differential operators with geometric
information. Their introduction also reports Oaku's calculation
\[
b_{I_2^{3,3}}(s)=(s+9/2)(s+4)(s+5),
\]
where the topological zeta function has poles only at $-9/2$ and $-4$.
This example already exhibits a Bernstein--Sato root beyond those
poles. For nonmaximal minors, one must determine the additional
roots and the repetitions of factors across the rank strata.

In our recent work of Chen–Zuo \cite{CZ25}, we develop the monodromy side of this
problem: We computed Verdier-monodromy zeta functions
along the exceptional divisor of the blow-up and proved the monodromy
conjecture for generic determinantal varieties. In that statement,
each relevant zeta pole $s_0$ yields a monodromy eigenvalue
$\exp(2\pi i s_0)$; the strong form asks that the pole itself be a
Bernstein--Sato root. Our later polynomial-PDE criterion of Chen and Zuo
\cite{CZPDE} detects the root $-N/d$ for a homogeneous ideal of
degree $d$ in $N$ variables under an explicit differential criterion,
and applies to maximal minors. Shi and Zuo \cite{SZ} give a
combinatorial root criterion and prove the pole-to-root form of the strong monodromy conjecture
for determinantal and Pfaffian ideals.

The analytic criterion of Chen and Zuo \cite{CZ} gives a longer interval
of roots: if the map defined by degree-$d$ generators in $N$
variables has differential of rank $N$ somewhere, then
\[
-\frac{N+u}{d},\qquad 0\le u<d,
\]
are roots. Applying this to nonmaximal minors and their transverse
rank slices produces the finite set in \eqref{eq:rootset} below.
Theorem~3.9 of \cite{CZ} proves the occurrence of these values, while
Remark~3.10 there leaves the full polynomial undetermined. We use this
root-existence theorem directly in Section~\ref{sec:existence}.
The new tasks are to exclude all further roots and to determine
the multiplicity of every root.

Multiplicity also has a broader interpretation through filtrations
on differential modules. For a holomorphic function,
L\H orincz and Yang \cite{LY} relate Bernstein--Sato-type
multiplicities to weight and Hodge filtrations and study the behavior
of the polynomials of powers.

For $r<m$, define the following set, with integer indices $i,u$:
\begin{equation}\label{eq:rootset}
\mathcal R=\left\{\frac{(m-i)(n-i)+u}{r-i}:
          0\le i<r,\quad 0\le u<r-i\right\}.
\end{equation}
For $\rho\in\C$, put
\[
J(\rho)=\{(i,u):0\le i<r,\ 0\le u<r-i,
                  \ (r-i)\rho=(m-i)(n-i)+u\}.
\]
Thus $\mathcal R$ is a set of distinct positive rational numbers,
whereas $J(\rho)$ retains all labels giving the same number. For a
polynomial $b$, $\mult_{s=-\rho}b$ denotes the exponent of $s+\rho$
in $b$, with value zero when $b(-\rho)\ne0$.

\begin{maintheorem}\label{thm:main}
If $1\le r<m\le n$, the roots of $b_{I_r}(s)$ are exactly $-\mathcal R$,
and
\begin{equation}\label{eq:multiplicity}
\mult_{s=-\rho}b_{I_r}(s)=
\begin{cases}
|J(\rho)|,&\rho\in\mathcal R\cap\Z,\\
1,&\rho\in\mathcal R\setminus\Z,\\
0,&\rho\notin\mathcal R.
\end{cases}
\end{equation}
Equivalently,
\begin{equation}\label{eq:mainproduct}
b_{I_r}(s)=
\prod_{i=0}^{r-1}\left(s+\left\lceil
              \frac{(m-i)(n-i)}{r-i}\right\rceil\right)
\prod_{\rho\in\mathcal R\setminus\Z}(s+\rho).
\end{equation}
The second product contains each distinct value once. For maximal minors,
\begin{equation}\label{eq:maximal}
b_{I_m^{m,n}}(s)=\prod_{j=n-m+1}^{n}(s+j).
\end{equation}
In particular, $b_{I_1^{m,n}}(s)=s+mn$.
\end{maintheorem}

Each fixed $i$ in \eqref{eq:rootset} gives $r-i$ consecutive numerator
values. Exactly one is divisible by $r-i$. Consequently, the first
product in \eqref{eq:mainproduct} counts the integral labels, including
their repetitions. For example, for $(m,n,r)=(6,8,4)$ the theorem gives
\begin{align*}
b_{I_4^{6,8}}(s)={}&(s+12)^3(s+15)
\left(s+\frac{35}{3}\right)\left(s+\frac{49}{4}\right)\\
&\hspace{12mm}\cdot\left(s+\frac{37}{3}\right)
\left(s+\frac{25}{2}\right)\left(s+\frac{51}{4}\right).
\end{align*}
The two labels for $25/2$ give a single factor, while the three labels
for $12$ give a triple factor.

Theorem~\ref{thm:main} completes the known root-existence results:
it proves that the displayed windows exhaust the roots and identifies
precisely which coincidences increase multiplicity. The distinction
between integral and nonintegral values is essential. Repeated
nonintegral labels still contribute only one factor; integral labels
contribute with their full multiplicity. This also explains why the
Bernstein--Sato polynomial contains more information than the set of
zeta poles.

Our second main result includes the multiplicity assertion in the
strong monodromy conjecture. For $0\le i<r$, put
\begin{equation}\label{eq:resolutionnumbers}
N_i=r-i,\qquad \nu_i=(m-i)(n-i),\qquad
p_r(\rho)=\#\{i:0\le i<r,\ \nu_i/N_i=\rho\}.
\end{equation}
Write $Z^{\mathrm{mot}}_{I_r}(T)$ for the motivic zeta function and
formally set $T=\mathbb L^{-s}$. The denominator convention for a
motivic pole at $s=-\rho$ and its order is specified in
Definition~\ref{def:motivicorder}. It avoids division in the
Grothendieck ring, which need not be an integral domain.

\begin{maintheorem}[Strong motivic monodromy with multiplicities]
\label{thm:strong}
For $1\le r\le m\le n$ and every positive rational number $\rho$,
\begin{equation}\label{eq:stronginequality}
\operatorname{ord}_{-\rho}Z^{\mathrm{mot}}_{I_r}
\le p_r(\rho)\le\mult_{s=-\rho}b_{I_r}(s).
\end{equation}
More precisely, the log-resolution formula for
$Z^{\mathrm{mot}}_{I_r}(T)$ has a common denominator
\begin{equation}\label{eq:motivicdenominator}
D_r(T)=\prod_{i=0}^{r-1}(1-\mathbb L^{-\nu_i}T^{N_i}),
\end{equation}
and the number of factors of any fixed slope $\nu_i/N_i=\rho$
in this denominator is at most $\mult_{s=-\rho}b_{I_r}(s)$.
Consequently, every motivic pole is a Bernstein--Sato root with at
least the same multiplicity. For the topological specialization,
\begin{equation}\label{eq:topologicalstrong}
b_{I_r}(s)Z^{\mathrm{top}}_{I_r}(s)\in\mathbb Q[s].
\end{equation}
\end{maintheorem}

Thus the result strengthens the previously established pole-to-root
statement by controlling orders, and does so at the level of the
motivic resolution formula. Cancellation of candidate poles can only
improve the inequality. The key arithmetic observation is that two
distinct indices in \eqref{eq:resolutionnumbers} can give the same
slope only when that slope is integral.

The integral multiplicities also have a sharp bound that depends
only on the size of the minors.
\begin{maintheorem}\label{thm:sharpbound}
For every positive integer $r$,
\begin{equation}\label{eq:sharpbound}
\max_{r\le m\le n}\ \max_{\rho\in\C}
      \mult_{s=-\rho}b_{I_r^{m,n}}(s)
=\left\lfloor\sqrt{4r-1}\right\rfloor.
\end{equation}
\end{maintheorem}
The upper bound is attained by an explicit family of rectangular
matrices; it is therefore a uniform optimal bound, not merely an
asymptotic estimate.

We have also verified the analogous root-and-multiplicity statement
for ideals of Pfaffians of a generic skew-symmetric matrix; that case
will be treated separately.

The proof of Theorem~\ref{thm:main} has three parts. Bordered minors
and the representation-theoretic description of ordinary powers
give relations in the graph module. Fourier transformation and a
Capelli identity then yield an annihilating polynomial. A weighted
Gaussian integral determines the total integral multiplicity, forcing
all integral bounds to be sharp. Finally, a filtration by kernels of
commuting operators has disjoint central supports at nonintegral
parameters, which forces the corresponding roots to be simple.

The paper is organized as follows. Section~\ref{sec:preliminaries}
recalls the graph definition of the ideal Bernstein--Sato polynomial,
weights and enveloping algebras, Schur decompositions and ordinary
powers, central characters and tensor products, and motivic zeta
functions in Sections~\ref{subsec:graph}--\ref{subsec:motivic}, respectively.
Section~\ref{sec:upper} proves the upper bound, through bordered-minor
identities, graph relations, the Fourier--Capelli calculation, and
the resulting annihilator in
Sections~\ref{subsec:border}--\ref{subsec:upperpolynomial}.
Section~\ref{sec:integer} determines the integral multiplicities:
Section~\ref{subsec:polebound} relates Mellin poles to Bernstein--Sato
multiplicities, Section~\ref{subsec:integerpole} constructs a pole of
order $r$, and Section~\ref{subsec:integermultiplicity} proves sharpness.
Section~\ref{sec:noninteger} constructs the bounded-word module in
Section~\ref{subsec:words} and separates its filtration layers in
Section~\ref{subsec:separation}, proving simplicity of nonintegral roots.
Section~\ref{sec:existence} quotes the root-existence results in
Section~\ref{subsec:knownroots} and completes the polynomial formula in
Section~\ref{subsec:completion}.
Finally, Section~\ref{sec:arithmetic} gives the interval count, proves
the sharp multiplicity bound, and establishes strong motivic monodromy
in Sections~\ref{subsec:interval}--\ref{subsec:strong}, respectively.


\begin{center}
\bf{Acknowledgements}
\end{center}

We are deeply grateful to Qiyue Tang for his indispensable contributions to this work through the use of artificial intelligence tools. His efforts played an essential role in the development of this paper. AI tools (ChatGPT and DeepSeek) were used during the preparation of the paper as technical assistants to help with checking calculations and arguments, developing details of proofs, and improving the exposition. All mathematical statements, proofs, and conclusions were verified by the authors. Y. Chen and H. Zuo were supported by BJNSF Grant 1252009. H. Zuo was supported by NSFC Grant 12671056.

\section{Preliminaries}\label{sec:preliminaries}
We collect the algebraic and representation-theoretic facts used in
the proof, together with the zeta-function conventions needed in the
last section.

\subsection{Bernstein--Sato polynomials and the graph module}\label{subsec:graph}
We use $\N=\{0,1,2,\ldots\}$. An empty product equals one, and
$[A,B]=AB-BA$ denotes the commutator of two operators.
Let $I=(f_1,\ldots,f_a)\subset S$ be a nonzero ideal. In the
determinantal application, $f_1,\ldots,f_a$ enumerate the $r$-minors,
so $a=\binom mr\binom nr$. Let $D_{X,t}$ be the Weyl algebra in the
variables $x_{ab}$ and $t_1,\ldots,t_a$, and let $D_X$ be the Weyl
algebra in the $x$ variables alone. A Weyl algebra is the algebra generated
by coordinates and their derivatives, subject to
$[\partial_{x_i},x_j]=\delta_{ij}$, with all coordinate-coordinate
and derivative-derivative commutators zero. For a multi-index $\alpha\in\N^a$,
write $|\alpha|=\sum_j\alpha_j$, $\alpha!=\prod_j\alpha_j!$,
$t^\alpha=\prod_jt_j^{\alpha_j}$, $f^\alpha=\prod_jf_j^{\alpha_j}$,
and $\partial_t^\alpha=\prod_j\partial_{t_j}^{\alpha_j}$.
Let $e_j$ be the $j$th standard multi-index.
The graph module of $(f_1,\ldots,f_a)$ is
\[
\mathcal M=S[\partial_{t_1},\ldots,\partial_{t_a}]\delta.
\]
Here $\delta$ is a formal generator, and for $h\in S$ and
$\alpha\in\N^a$, the action is
\begin{align}
\partial_{x_\ell}(h\partial_t^\alpha\delta)
 &= (\partial_{x_\ell}h)\partial_t^\alpha\delta
       -\sum_j h(\partial_{x_\ell}f_j)
                       \partial_t^{\alpha+e_j}\delta,\label{eq:graphdx}\\
t_j(h\partial_t^\alpha\delta)
 &=hf_j\partial_t^\alpha\delta
       -\alpha_jh\partial_t^{\alpha-e_j}\delta.\label{eq:grapht}
\end{align}
The index $\ell$ in \eqref{eq:graphdx} ranges over the matrix entries;
all sums over generator indices range from $1$ to $a$.
A term with $\alpha_j=0$ in \eqref{eq:grapht} is zero.

Give $t_j$ degree $1$, $\partial_{t_j}$ degree $-1$, and the $x$-Weyl
operators degree zero. The decreasing filtration $V^\ell D_{X,t}$
is spanned by normally ordered monomials of degree at least $\ell$,
for $\ell\in\Z$.
Set
\[
U^\ell=V^\ell D_{X,t}\delta,\qquad Q=U^0/U^1,\qquad
s=-\sum_j\partial_{t_j}t_j.
\]
The module $Q$ is cyclic over $\gr_V^0D_{X,t}=V^0D_{X,t}/V^1D_{X,t}$,
with generator the class of $\delta$. The operator $s$ is central in
this ring. By \cite[Sections~2.4--2.5]{BMS}, there is a nonzero polynomial
annihilating the action of $s$ on $Q$. The monic polynomial of least
degree with this property is the \emph{unshifted ideal Bernstein--Sato
polynomial} $b_I(s)$.
In particular, a polynomial in $s$ kills $Q$ if and only if it kills
the class of $\delta$. On the free graph module,
\begin{equation}\label{eq:saction}
s(h\partial_t^\alpha\delta)=
-\sum_jhf_j\partial_t^{\alpha+e_j}\delta
       +|\alpha|h\partial_t^\alpha\delta.
\end{equation}
Let $E=\sum_{a,b}x_{ab}\partial_{ab}$ be the Euler operator in $X$,
where $\partial_{ab}=\partial/\partial x_{ab}$. For the $r$-minor generators, homogeneity gives
$E\delta=rs\delta$.

For a principal ideal $(f)$, this definition agrees with the
classical polynomial $b_f(s)$. In that setting $f^s$ is the formal
generator of $S[f^{-1},s]f^s$, with
\[
\partial_{x_i}(h f^s)=
\left(\partial_{x_i}h+s h\frac{\partial_{x_i}f}{f}\right)f^s.
\]
The equation in the introduction is interpreted in this module.
The ideal polynomial is independent of its generators, is unchanged
by a smooth factor, and is the least common multiple of the local
polynomials; see \cite[Sections~2.4--2.5]{BMS}. For an equidimensional
closed subvariety $Y$ of codimension $c$, the embedding-independent
normalization is $b_Y(s)=b_{I_Y}(s-c)$; we use $b_I$, without this
shift, throughout our main theorems.

We will also use the comparison with a single polynomial. If
$g(x,y)=\sum_{j=1}^a y_j f_j(x)$, then
\begin{equation}\label{eq:mustata}
b_g(s)=(s+1)b_I(s)
\end{equation}
by \cite[Theorem~1.1]{Mustata}.

\subsection{Enveloping algebras and weights}\label{subsec:weights}
Let $\g=\mathfrak{gl}_m(\C)$, with matrix-unit basis $e_{uv}$ and
bracket
\[
[e_{uv},e_{ij}]=\delta_{vi}e_{uj}-\delta_{ju}e_{iv}.
\]
Its universal enveloping algebra is the associative algebra
\[
U(\g)=T(\g)/\langle a\otimes b-b\otimes a-[a,b]:a,b\in\g\rangle,
\]
where $T(\g)=\bigoplus_{k\ge0}\g^{\otimes k}$ is the tensor algebra.
A representation of $\g$ extends uniquely to a unital representation
of $U(\g)$. The center is
$Z=Z(U(\g))=\{z\in U(\g):zu=uz\text{ for every }u\in U(\g)\}$.

Let $\mathfrak h$ be the diagonal subalgebra. A vector $w\ne0$ has
weight $\lambda=(\lambda_1,\ldots,\lambda_m)\in\C^m$ if
$e_{uu}w=\lambda_u w$ for all $u$. It is a highest-weight vector
relative to the upper triangular matrices if also $e_{uv}w=0$ for
$u<v$. A highest-weight module of highest weight $\lambda$ is a
module generated by such a vector. A finite-dimensional algebraic
$GL_m$-representation is a homomorphism of algebraic groups
$GL_m\to GL(V)$; differentiating it gives a $\g$-module. Its weights
are integral, and it is a direct sum of irreducible representations.
An irreducible polynomial representation is uniquely determined
by a highest weight $\lambda_1\ge\cdots\ge\lambda_m\ge0$.

On polynomials in the matrix entries, the representation used below is
\begin{equation}\label{eq:rowaction}
e_{uv}\longmapsto A_{uv}=\sum_{b=1}^n x_{ub}\partial_{vb}.
\end{equation}
Indeed, $[A_{uv},A_{ij}]=\delta_{vi}A_{uj}-\delta_{ju}A_{iv}$.
The notation means that the abstract element $e_{uv}$ acts by the
differential operator $A_{uv}$. In particular,
$\sum_u A_{uu}=E$ is the image of the trace element of $\g$.
On a tensor product the Lie algebra acts diagonally:
$a(v\otimes w)=(av)\otimes w+v\otimes(aw)$.
These conventions agree with the representation-theoretic setup of
\cite[Sections~2.1--2.3]{LRWW}.

\subsection{Schur decompositions and ordinary powers}\label{subsec:schur}
We put $I_0=S$ and $I_t=0$ for $t>m$. A partition is a nonincreasing finite sequence of
nonnegative integers, extended by zeros; its size is
$|\lambda|=\sum_i\lambda_i$ and its length $\ell(\lambda)$ is the number
of positive parts;
$\lambda_1$ is its largest part. The Schur representation
$S_\lambda\C^m$ is the irreducible polynomial $GL_m$-representation
of highest weight $\lambda$; representations indexed by partitions whose length
exceeds the dimension are zero. Choose the row and column actions
so that the linear variables transform as $\C^m\otimes\C^n$.
The Cauchy decomposition gives, for each $d\ge0$,
\begin{equation}\label{eq:cauchy}
S_d=\operatorname{Sym}^d(\C^m\otimes\C^n)
=\bigoplus_{\substack{|\lambda|=d\\\ell(\lambda)\le m}}
S_\lambda\C^m\otimes S_\lambda\C^n;
\end{equation}
each summand occurs once \cite[(2.1)]{Raicu}.
Here $S_d$ denotes the homogeneous degree-$d$ part of $S$ and
$\operatorname{Sym}^d$ denotes the $d$th symmetric power.
The span of the $t$-minors is
$\bigwedge^t\C^m\otimes\bigwedge^t\C^n$; the weights of
$\bigwedge^t\C^m$ have $t$ entries equal to one and all others zero.

For $v\in\N$, $I_r^v$ means the \emph{ordinary} ideal power, generated
by products of $v$ minors. Its degree-$rv$ component is
\cite[(5.4)]{Raicu}
\begin{equation}\label{eq:ordinary}
(I_r^v)_{rv}=
\bigoplus_{\substack{|\lambda|=rv\\\lambda_1\le v}}
S_\lambda\C^m\otimes S_\lambda\C^n.
\end{equation}

\subsection{Capelli operators and central characters}\label{subsec:central}
An algebra homomorphism $\chi:Z\to\C$ is a \emph{central character}.
A module has exact central character $\chi$ if $z$ acts as the scalar
$\chi(z)$ for every $z\in Z$. Every highest-weight module of highest
weight $\lambda$ has such a character, denoted $\chi_\lambda$.

For $0\le p\le m$, define the Capelli operators
\begin{equation}\label{eq:capellidefinition}
C_0=1,\qquad
C_p=\sum_{\substack{I\subseteq\{1,\ldots,m\},\ |I|=p\\
                   J\subseteq\{1,\ldots,n\},\ |J|=p}}
       \Delta_{I,J}(X)\Delta_{I,J}(\partial),
\end{equation}
where $\Delta_{I,J}$ is the determinant of the indicated submatrix,
with row and column indices in increasing order. Multiplication by
the first determinant precedes the constant-coefficient derivative
operator in this formula. In particular, $C_1=E$.
There are algebraically independent generators of $Z$ whose images
under \eqref{eq:rowaction} are $C_1,\ldots,C_m$; we use the same
symbols for these central elements and their images.
On a highest-weight module of highest weight
$\lambda=(\lambda_1,\ldots,\lambda_m)$, their values are
\begin{equation}\label{eq:capellieigenvalue}
c_p(\lambda)=\sum_{1\le i_1<\cdots<i_p\le m}
                         \prod_{a=1}^p(\lambda_{i_a}+p-a).
\end{equation}
The central character $\chi_\lambda:Z\to\C$ is determined by the
unordered shifted coordinates
\begin{equation}\label{eq:HC}
\{\lambda_1+m-1,\lambda_2+m-2,\ldots,\lambda_m\}.
\end{equation}
These are the Capelli and Harish--Chandra conventions of
\cite[(0.9), (0.10), Example~6.2, Proposition~6.5, Theorem~6.7]{OO}.
The highest homogeneous terms of \eqref{eq:capellieigenvalue} are the
elementary symmetric polynomials; filtered induction gives the
generation of $Z$. The formulas extend from dominant integral weights
to all complex weights by polynomiality and Zariski density.

For a $Z$-module $M$, write
$\Ann_Z M=\{z\in Z:zw=0\text{ for every }w\in M\}$.
If it is annihilated by an ideal of finite codimension,
its central support is the finite set of maximal ideals of
$Z/\Ann_Z M$, or equivalently the central characters occurring in
its generalized-character decomposition. The summand for a character
$\chi:Z\to\C$ consists of the vectors killed by a power of
$\ker\chi$. The decomposition follows from the Chinese remainder
theorem for the finite-dimensional algebra $Z/\Ann_Z M$.

\begin{lemma}[{\cite[Theorem~2.5, Corollary~2.6, and p.~260]{BG}}]\label{lem:tensorsupport}
Let $M$ be a $U(\g)$-module annihilated by a finite-codimensional
ideal of $Z$, and let $V$ be a finite-dimensional algebraic
$GL_m$-representation. Then $V\otimes M$ is annihilated by a
finite-codimensional ideal of $Z$. If $M$ has central support
$\{\chi_\lambda\}$, the support of $V\otimes M$ is contained in
\[
\{\chi_{\lambda+\nu}:\nu\text{ is a weight of }V\}.
\]
\end{lemma}
\begin{proof}
For a semisimple Lie algebra, let $Z_{\mathrm{ss}}$ be the center
of its universal enveloping algebra
and let $\mathcal I_V\subset Z_{\mathrm{ss}}\otimes Z_{\mathrm{ss}}$
be the ideal of relations between the output and input central
actions of $V\otimes M$. Theorem~2.5 of \cite{BG} identifies its
zero set with the weight-shift correspondences and proves that
these relations annihilate $V\otimes M$ for every module $M$.
Corollary~2.6 shows that the quotient by $\mathcal I_V$ is finite
over the input center. If a finite-codimensional ideal $J$ kills
$M$, the further quotient by $Z_{\mathrm{ss}}\otimes J$ is finite
dimensional. The kernel of the output-center action on this
quotient is a finite-codimensional ideal killing $V\otimes M$.
Its support lies over the input characters, with output coordinates
shifted by weights of $V$, as specified on p.~260 of \cite{BG}.

For $m=1$, decompose $V$ into its trace-weight spaces. For $m>1$,
pass from $\mathfrak{sl}_m$ to $\mathfrak{gl}_m$ by retaining the
trace element $E=\sum_u A_{uu}$. Write ${\bf1}=(1,\ldots,1)\in\C^m$.
If $\tau=\sum_i\lambda_i$ and a summand of
$V$ has trace weight $j$, the $\mathfrak{sl}_m$ parameter is
\[
\lambda-\frac{\tau}{m}{\bf1}
+(m-1,m-2,\ldots,0)-\frac{m-1}{2}{\bf1}.
\]
For a weight $\nu$ of that summand, $\sum_i\nu_i=j$.
Adding its traceless part and restoring the trace produces precisely
the shifted coordinates of $\lambda+\nu$. Moreover, on that tensor
summand, for every integer $\ell\ge0$,
\[
(E-(\tau+j))^\ell=1\otimes(E-\tau)^\ell.
\]
Thus a uniform annihilating power for the trace is retained, as is
the finite-codimensional annihilator supplied by the semisimple
result. The finite-support version follows by decomposing $M$ into
its generalized-character summands.
\end{proof}

\subsection{Motivic zeta functions and pole orders}\label{subsec:motivic}
The Grothendieck ring $K_0(\operatorname{Var}_{\C})$ is generated by
isomorphism classes $[Y]$ of complex algebraic varieties, with
$[Y]=[W]+[Y\setminus W]$ for a closed subvariety $W\subseteq Y$
and $[Y][Y']=[Y\times Y']$. Put
$\mathbb L=[\mathbb A^1]$ and
$\mathcal M_{\C}=K_0(\operatorname{Var}_{\C})[\mathbb L^{-1}]$.
For $X=\mathbb A^d$, an arc is a tuple in $\C[[t]]^d$; its
$q$-jet is its reduction modulo $t^{q+1}$. Let
\[
\mathcal X_{q,I}=\{\gamma\in J_q(X):
           \min_j\operatorname{ord}_t f_j(\gamma)=q\},
\]
where $J_q(X)=\mathbb A^{d(q+1)}$. We use the normalized,
non-equivariant motivic zeta function
\begin{equation}\label{eq:motivicdefinition}
Z_I^{\mathrm{mot}}(T)
=\sum_{q\ge0}[\mathcal X_{q,I}]\mathbb L^{-d(q+1)}T^q
\in\mathcal M_{\C}[[T]].
\end{equation}
The order of the zero truncated series is understood to be $>q$.
Multiplying this convention by $\mathbb L^d$, or omitting its
constant term, does not affect the denominator bounds below.

A log resolution of $I$ is a proper birational morphism
$\pi:Y\to X$ with $Y$ smooth such that
$I\mathcal O_Y=\mathcal O_Y(-\sum_{i\in A}N_iE_i)$ and the relevant
divisors have simple normal crossings: locally their union is a
union of coordinate hyperplanes. Write
$K_{Y/X}=\sum_{i\in A}(\nu_i-1)E_i$, and
\[
E_J^\circ=\left(\bigcap_{i\in J}E_i\right)
                \setminus\bigcup_{i\notin J}E_i,
\qquad E_\varnothing^\circ=Y\setminus\bigcup_{i\in A}E_i.
\]
For a principalization that is an isomorphism off $V(I)$, all the
relevant $N_i$ are positive. The change-of-variables formula for
motivic integration gives
\begin{equation}\label{eq:motivicresolution}
Z_I^{\mathrm{mot}}(T)=\mathbb L^{-d}
 \sum_{J\subseteq A}[E_J^\circ]
 \prod_{i\in J}
 \frac{(\mathbb L-1)\mathbb L^{-\nu_i}T^{N_i}}
      {1-\mathbb L^{-\nu_i}T^{N_i}}.
\end{equation}
We use this standard resolution formula in the ideal version
recorded in \cite[Theorem~2.7]{CZ25}; for the underlying motivic
change of variables see \cite[Lemma~3.3]{DL}. The factor $\mathbb L^{-d}$
corresponds to the normalization in \eqref{eq:motivicdefinition}.
The topological zeta function is
\begin{equation}\label{eq:topologicaldefinition}
Z_I^{\mathrm{top}}(s)=
\sum_{J\subseteq A}\chi_c(E_J^\circ)
                    \prod_{i\in J}\frac1{N_i s+\nu_i},
\end{equation}
where $\chi_c$ is the compactly supported Euler characteristic.
Both formulas are independent of the chosen resolution.

\begin{definition}[Denominator convention for motivic pole orders]
\label{def:motivicorder}
Consider a rational motivic series expressible as a finite sum
\begin{equation}\label{eq:denominatorpresentation}
Z(T)=\sum_{\alpha}P_\alpha(T)
          \prod_{j=1}^{h_\alpha}
           (1-\mathbb L^{-a_{\alpha j}}T^{b_{\alpha j}})^{-1},
\qquad
P_\alpha(T)\in\mathcal M_{\C}[T],\quad
 a_{\alpha j},b_{\alpha j}\in\Z_{>0}.
\end{equation}
For a positive rational number $\rho$, its pole order at
$s=-\rho$ is the least nonnegative integer $e$ for which such a
presentation exists with at most $e$ factors satisfying
$a_{\alpha j}/b_{\alpha j}=\rho$ in each summand, repetitions counted.
A positive order means that $-\rho$ is a pole. This is a denominator
notion of order in the motivic ring. After Hodge realization, with $\mathbb L$ sent to the indeterminate
$uv$, the ordinary order along $T=(uv)^\rho$ (in a suitable finite
field extension) is bounded by this order.
\end{definition}

This convention states an upper bound without assuming unique
factorization in $\mathcal M_{\C}$. Formula~\eqref{eq:motivicresolution}
shows that the order at $-\rho$ is at most the number of divisors
with $\nu_i/N_i=\rho$. Theorem~\ref{thm:strong} will prove the
stronger assertion that one common denominator works simultaneously
with the required bounds at every slope. For
\eqref{eq:topologicaldefinition}, pole orders have their usual meaning
for a rational function of the complex variable $s$.

\section{An upper bound for the Bernstein--Sato polynomial}\label{sec:upper}
From now on, the generators $f_j$ are the $r$-minors. We first
construct an annihilator of their graph quotient; this bounds both
the possible roots and their multiplicities.

\subsection{Bordered-minor identities}\label{subsec:border}

For an entry index $\eta=(a,b)$ define
\[
R_{ab}=\sum_{i=1}^m\sum_{j=1}^n x_{ib}x_{aj}\partial_{ij},
\qquad F_{ab}=R_{ab}-sx_{ab}.
\]
For ordered row and column lists $I,J$ of the same length, write
$\Delta_{I,J}(X)=\det X_{I,J}$ and define the bordering operation by
\[
\mathfrak b_{ab}\Delta_{I,J}
=\det\begin{pmatrix}X_{I,J}&X_{I,b}\\ X_{a,J}&x_{ab}\end{pmatrix}.
\]
A determinant with a repeated row or column is zero. Repeated bordering
appends the specified rows and columns in the same order. Swapping two
bordering operations swaps both rows and columns, so these operations
commute.

\begin{lemma}\label{lem:bell}
The operators $F_{ab}$ commute. If $\eta_1,\ldots,\eta_k$ are entry
indices, then
\begin{equation}\label{eq:bell}
F_{\eta_1}\cdots F_{\eta_k}\delta
=\sum_\pi(-1)^{k-|\pi|}
 \sum_{(j_D)_{D\in\pi}}
 \left(\prod_{D\in\pi}\mathfrak b_D f_{j_D}\right)
 \left(\prod_{D\in\pi}\partial_{t_{j_D}}\right)\delta.
\end{equation}
Here $\pi$ runs over the set partitions of $\{1,\ldots,k\}$,
each $j_D$ runs from $1$ to $a$, and
$\mathfrak b_D=\prod_{\ell\in D}\mathfrak b_{\eta_\ell}$.
For $k=0$, the empty partition gives $\delta$.
\end{lemma}
\begin{proof}
For a square submatrix $M=X_{I,J}$, set $u=X_{I,b}$ and
$w=X_{a,J}$. The determinant identity
\[
\det\begin{pmatrix}M&u\\w&x_{ab}\end{pmatrix}
=x_{ab}\det M-w\operatorname{adj}(M)u
\]
holds polynomially, including singular $M$; $\operatorname{adj}(M)$
is the transpose of the cofactor matrix. On $M$, the vector field
$R_{ab}$ has variation $uw$. Differentiating its determinant gives
$R_{ab}\det M=w\operatorname{adj}(M)u$. Thus, for every minor $\Delta$,
\begin{equation}\label{eq:border}
R_{ab}\Delta=x_{ab}\Delta-\mathfrak b_{ab}\Delta.
\end{equation}
On each coordinate $x_{ij}$,
\[
R_{ab}R_{cd}x_{ij}
=x_{ib}x_{ad}x_{cj}+x_{id}x_{cb}x_{aj}
=R_{cd}R_{ab}x_{ij}.
\]
The commutator of two vector fields is a derivation, so its vanishing
on all coordinates gives $[R_{ab},R_{cd}]=0$. Also
$R_{ab}x_{cd}=x_{cb}x_{ad}=R_{cd}x_{ab}$.
Since $s$ commutes with $x$ and $\partial_x$,
\[
[F_{ab},F_{cd}]
=[R_{ab},R_{cd}]-s(R_{ab}x_{cd}-R_{cd}x_{ab})=0.
\]

Equations~\eqref{eq:saction} and \eqref{eq:border} give
\[
F_\eta(h\partial_t^\alpha\delta)
=(R_\eta h-|\alpha|x_\eta h)\partial_t^\alpha\delta
 +\sum_jh(\mathfrak b_\eta f_j)\partial_t^{\alpha+e_j}\delta.
\]
For $h=1$ and $\alpha=0$, this gives
$F_\eta\delta=\sum_j(\mathfrak b_\eta f_j)\partial_{t_j}\delta$.
Suppose a term on the right of \eqref{eq:bell} has $v$ blocks. Its
coefficient is a product of $v$ minors, and its $\partial_t$-degree is
$v$. Applying $R_\eta$ to that product gives $v x_\eta$ times the
product, together with the negative of each possible additional
border. The former term cancels the $-v x_\eta$ contribution from
$-sx_\eta$. The graph-derivative terms give one additional singleton
block. Adding a border to an existing block changes the sign, whereas
adding a singleton does not. This is the recursion for the set
partitions in \eqref{eq:bell}, and proves the formula by induction.
\end{proof}

\subsection{Relations in the graph quotient}\label{subsec:graphrelations}

\begin{proposition}\label{prop:graphrelations}
For $0\le k\le r$ and all entry indices $\eta_1,\ldots,\eta_k$,
\begin{equation}\label{eq:graphrelations}
I_{r-k}(X)F_{\eta_1}\cdots F_{\eta_k}\delta=0
\quad\text{in }Q.
\end{equation}
In particular, every product of $r$ of the $F_{ab}$ kills $\delta$ in $Q$.
\end{proposition}
\begin{proof}
For $k=0$, a generator $f_j\delta=t_j\delta$ belongs to $U^1$.
For $k>0$, consider a term with $v$ blocks in \eqref{eq:bell}. The
sizes of its minors are $r+j_1,\ldots,r+j_v$, where each $j_\ell\ge1$
and $\sum j_\ell=k$. Multiplication by an $(r-k)$-minor gives total
degree $r(v+1)$ and at most $v+1$ exterior-power factors on each
matrix side. Each weight coordinate in an exterior power is zero or
one. Hence every highest weight in their tensor product has largest
part at most $v+1$. Equation~\eqref{eq:ordinary} implies
\[
I_{r-k}I_{r+j_1}\cdots I_{r+j_v}\subseteq I_r^{v+1}.
\]
For $k=r$, the factor $I_0$ contributes no exterior-power factor,
so the same bound holds with only $v$ factors.

Write a coefficient in $I_r^{v+1}$ as a sum of polynomial multiples
of $f^\beta$, with $|\beta|=v+1$. For $|\alpha|=v$, the graph action gives
\[
f^\beta\partial_t^\alpha\delta
=\partial_t^\alpha t^\beta\delta\in U^1,
\]
since $|\beta|-|\alpha|=1$. All the terms in \eqref{eq:bell} therefore
vanish after multiplication by $I_{r-k}$ and passage to $Q$.
\end{proof}

\subsection{Fourier transformation and a Capelli identity}\label{subsec:capelli}

Apply the Weyl algebra automorphism
$x_{ab}\mapsto\partial_{ab}$, $\partial_{ab}\mapsto-x_{ab}$, $s\mapsto s$
to the relations of the $D_X[s]$-module generated by $\delta$ in $Q$.
Denote the transformed generator by $v$. Introduce the polynomials
\[
q=-rs-mn,\qquad z=s+m+n,
\]
which are the transformed Euler eigenvalue and the coefficient in the
following operators:
\begin{equation}\label{eq:B}
B_{ab}=\sum_{i,j}x_{ij}\partial_{ib}\partial_{aj}+z\partial_{ab}.
\end{equation}
Indeed, the transform of $F_{ab}$ is $-B_{ab}$, because
\[
\sum_{i,j}\partial_{ib}\partial_{aj}x_{ij}
=\sum_{i,j}x_{ij}\partial_{ib}\partial_{aj}+(m+n)\partial_{ab}.
\]
The $B_{ab}$ commute, and $[E,B_{ab}]=-B_{ab}$. Here and below $B^k v$
means any word of length $k$ in these operators applied to $v$.
Write $I_t(\partial)$ for the constant-coefficient operators
obtained by substituting $\partial_{ab}$ for $x_{ab}$ in $I_t$.
For $0\le k\le r$, Proposition~\ref{prop:graphrelations} becomes
\begin{equation}\label{eq:fourierrelations}
Ev=qv,\qquad EB^kv=(q-k)B^kv,\qquad
I_{r-k}(\partial)B^kv=0,\qquad B^rv=0.
\end{equation}

We use the Capelli operators of \eqref{eq:capellidefinition}.

\begin{lemma}\label{lem:capelli}
For $0\le p<m$, one has the identity of differential operators
\begin{equation}\label{eq:capelli}
\sum_{a,b}[C_p,x_{ab}]B_{ab}
    =C_p(E+pz-p^2)-(p+1)C_{p+1}.
\end{equation}
Consequently, if $w$ belongs to a $D_X[s]$-module and
$B_{ab}w=0$ for every $a,b$, then
\begin{equation}\label{eq:commonkernel}
t!C_tw=\prod_{i=0}^{t-1}(E+iz-i^2)w\qquad(0\le t\le m).
\end{equation}
\end{lemma}
\begin{proof}
Use commuting variables $\Xi=(\xi_{ab})$ to denote normal symbols.
Write polynomial coefficients to the left of all derivatives and
replace each $\partial_{ab}$ by $\xi_{ab}$ to obtain the normal symbol.
If two operators have normal symbols $A(X,\Xi)$ and $B(X,\Xi)$,
their product has symbol
\[
\sum_{\alpha\in\N^{mn}}\frac1{\alpha!}
 (\partial_\Xi^\alpha A)(\partial_X^\alpha B),
\qquad \alpha!=\prod_{a,b}\alpha_{ab}!.
\]
Here $\delta_{ij}$ is one if $i=j$ and zero otherwise. The formula
follows by moving derivatives past coefficients using
$[\partial_{ab},x_{ij}]=\delta_{ai}\delta_{bj}$.
For a square matrix $A$, let $e_p(A)$ be the sum of its principal
$p$-minors, with $e_0(A)=1$. A superscript $t$ on a matrix denotes
transpose. Put
\[
P_p=e_p(X\Xi^t)=\sum_{|I|=|J|=p}\Delta_{I,J}(X)\Delta_{I,J}(\Xi),
\quad e=\sum_{a,b}x_{ab}\xi_{ab}.
\]
Also put
$K_{ab}=\sum_{i,j}x_{ij}\xi_{ib}\xi_{aj}$. The symbol of
$[C_p,x_{ab}]$ is $\partial_{\xi_{ab}}P_p$.

The variation $\dot\Xi=\Xi X^t\Xi$ gives
$\frac{d}{dt}(X\Xi^t)=(X\Xi^t)^2$. Newton's identity thus yields
\begin{equation}\label{eq:newtonsymbol}
\sum_{a,b}(\partial_{\xi_{ab}}P_p)K_{ab}
          =eP_p-(p+1)P_{p+1}.
\end{equation}
For instance, this follows first for diagonal matrices from
$\sum_j a_j^2 e_{p-1}(a_1,\ldots,\widehat a_j,\ldots,a_m)
=e_1e_p-(p+1)e_{p+1}$, and then as a polynomial identity.

Because $K_{ab}$ is linear in $X$, normal composition contributes
exactly one further term:
\begin{equation}\label{eq:hessian}
\sum_{a,b,i,j}
 (\partial_{\xi_{ij}}\partial_{\xi_{ab}}P_p)\xi_{ib}\xi_{aj}
                       =-p(p-1)P_p.
\end{equation}
To verify it, expand each determinant $\Delta_{I,J}(\Xi)$. The two
derivatives choose two distinct slots. Replacing the selected entries
by $\xi_{ib}\xi_{aj}$ exchanges their columns and changes the
determinant sign. There are $p(p-1)$ ordered pairs of slots; repeated
rows or columns contribute zero. Summing proves \eqref{eq:hessian}.
The term $z\partial_{ab}$ contributes $pzP_p$ by homogeneity in $\Xi$.
Finally, the normal symbol of $C_pE$ is $eP_p+pP_p$. These identities
give \eqref{eq:capelli}, including $p=0$.

If every $B_{ab}$ kills $w$, the left side of \eqref{eq:capelli}
kills $w$. Each $C_p$ commutes with $E$ and $z$. Induction from
$C_0=1$ therefore gives \eqref{eq:commonkernel}.
\end{proof}

\subsection{The annihilating polynomial}\label{subsec:upperpolynomial}

\begin{proposition}\label{prop:upper}
For $1\le r\le m\le n$,
\begin{equation}\label{eq:upper}
b_{I_r}(s)\ \bigm|\
\prod_{i=0}^{r-1}\prod_{u=0}^{r-i-1}
       \big((r-i)s+(m-i)(n-i)+u\big).
\end{equation}
\end{proposition}
\begin{proof}
For $0\le k<r$ set
$P_k(q,z)=\prod_{i=0}^{r-k-1}(q-k+iz-i^2)$.
Starting with $B^rv=0$, descend through $k=r-1,\ldots,0$.
Suppose that $w$ is a polynomial in $s$ times $v$, and that every
word of length $k+1$ kills $w$. For any word of length $k$, put
$u=B^kw$. Equations~\eqref{eq:fourierrelations} give
\[
B_{ab}u=0\quad\text{for all }a,b,\qquad
C_{r-k}u=0,\qquad Eu=(q-k)u.
\]
Equation~\eqref{eq:commonkernel} gives $P_k(q,z)u=0$.
Replacing $w$ by $P_k(q,z)w$ therefore makes every word of length $k$
vanish. At the end of the descent,
\[
\prod_{k=0}^{r-1}P_k(q,z)v=0.
\]
The factors satisfy
\[
q-k+iz-i^2=-\big((r-i)s+(m-i)(n-i)+k\big).
\]
The pairs $0\le k<r$, $0\le i<r-k$ are the pairs
$0\le i<r$, $0\le k<r-i$. Transform back to $Q$ and use its cyclicity
and the centrality of $s$ to obtain \eqref{eq:upper}.
\end{proof}

\section{Multiplicities of the integral roots}\label{sec:integer}

Write $\mu(\rho)=\mult_{s=-\rho}b_{I_r}(s)$. We first relate a pole
at a negative integer to the total of these multiplicities over $\Z$.
All integrals use ordinary Lebesgue measure on the relevant complex
vector space, so that $\int_{\C}e^{-|z|^2}\,dz=\pi$.
A Schwartz function is smooth, and every derivative decays faster
than any inverse polynomial. We use the Gamma function
$\Gamma(z)=\int_0^\infty t^{z-1}e^{-t}\,dt$ for $\Re z>0$ and its
meromorphic continuation, determined by $\Gamma(z+1)=z\Gamma(z)$.

\subsection{From Mellin poles to integral multiplicities}\label{subsec:polebound}

\begin{lemma}\label{lem:polebound}
Let $I=(f_1,\ldots,f_a)\subset\C[z_1,\ldots,z_N]$ be nonzero,
and let $\phi$ be a Schwartz function on $\C^N$.
Define, initially for $\Re s>0$,
\[
Z_I(\phi;s)=\int_{\C^N}\left(\sum_{j=1}^a|f_j(z)|^2\right)^s
                         \phi(z)\,dz.
\]
If its meromorphic continuation has a pole of order $R$ at a negative
integer $-\rho$, then
\begin{equation}\label{eq:totalpolebound}
R\le\sum_{u\in\Z}\mult_{s=-u}b_I(s).
\end{equation}
\end{lemma}
\begin{proof}
For a nonzero polynomial $g$ and a Schwartz test function $\psi$
on its ambient complex vector space, choose a Bernstein equation
$P(s)g^{s+1}=b_g(s)g^s$. For a positive integer $N$, put
$P_N(s)=P(s)P(s+1)\cdots P(s+N-1)$.
Iterating the equation and integrating by parts in the holomorphic
variables gives
\begin{equation}\label{eq:onesided}
\left(\prod_{j=0}^{N-1}b_g(s+j)\right)
       \int |g|^{2s}\psi
=\int g^N|g|^{2s}P_N(s)^{\mathsf t}\psi,
\end{equation}
The transpose $\mathsf t$ is the bilinear formal transpose:
$\partial_{z_i}^{\mathsf t}=-\partial_{z_i}$, with multiplication
and the scalar parameter $s$ fixed. One may first use
cutoffs and $\Re s$ large, where differentiation across the zero set
is justified, and then continue meromorphically. The right side is
holomorphic whenever $N+2\Re s>0$. Near $g=0$, positive powers of
$|g|$ control every parameter derivative involving $\log|g|$; at
infinity the Schwartz decay controls the polynomial coefficients and
growth of $g$. Thus the pole order on the left is bounded by the
zero order of the single shifted product.

Apply this to
$g(z,y)=\sum_{j=1}^a y_jf_j(z)$ and
$\psi(z,y)=\phi(z)e^{-\|y\|^2}$, where $y\in\C^a$ and
$\|y\|^2=\sum_j|y_j|^2$.
The comparison theorem \eqref{eq:mustata} gives
$b_g(s)=(s+1)b_I(s)$.
For fixed $z$ outside the common zero set of the generators, a unitary
change in $y$ makes $g$ equal to $(\sum_j|f_j(z)|^2)^{1/2}y_1$.
The other $a-1$ coordinates contribute $\pi^{a-1}$, while
polar coordinates give
\[
\int_{\C}|y_1|^{2s}e^{-|y_1|^2}\,dy_1
=2\pi\int_0^\infty t^{2s+1}e^{-t^2}\,dt=\pi\Gamma(s+1).
\]
The common zero set has measure zero. Fubini's theorem for
$\Re s>0$ and meromorphic continuation therefore give
\begin{equation}\label{eq:gaussiancompare}
\int |g|^{2s}\psi=\pi^a\Gamma(s+1)Z_I(\phi;s).
\end{equation}
At $s=-\rho$, the left side has order $R+1$. Choose an integer
$N>2\rho$. Among the extra factors $s+j+1$, exactly one vanishes
there. Equations~\eqref{eq:onesided}--\eqref{eq:gaussiancompare} imply
\[
R+1\le 1+\sum_{j=0}^{N-1}\mult_{s=-(\rho-j)}b_I(s)
\le1+\sum_{u\in\Z}\mult_{s=-u}b_I(s).
\]
This proves \eqref{eq:totalpolebound}.
\end{proof}

\subsection{A weighted Gaussian pole of order \texorpdfstring{$r$}{r}}\label{subsec:integerpole}

For variables $w_1,\ldots,w_t$, write
$\Delta_t(w)=\prod_{i<j}(w_i-w_j)$ and let $e_j(w)$ be the $j$th
elementary symmetric polynomial. An empty product is one. For a
polynomial or power series $H(w)$, the notation
$[w_1^{a_1}\cdots w_t^{a_t}]H$ denotes the coefficient of that
monomial, with value zero if the monomial does not occur.

\begin{proposition}\label{prop:integerpole}
For $2\le r<m\le n$, let $X^*$ be the conjugate transpose and put
$\|X\|^2=\sum_{a,b}|x_{ab}|^2$. Define
\[
\phi(X)=\det(XX^*)^r e^{-\|X\|^2},\qquad
\rho=(m-r+1)(n+1).
\]
Then $Z_{I_r}(\phi;s)$ has a pole
of order exactly $r$ at $s=-\rho$.
\end{proposition}
\begin{proof}
The function $\phi$ is a polynomial in the real coordinates times a
Gaussian, and hence is Schwartz. If
$\lambda_1>\cdots>\lambda_m>0$ are the eigenvalues of $XX^*$,
Cauchy--Binet gives $\sum_j|f_j(X)|^2=e_r(\lambda)$.
Write $Z(s)=Z_{I_r}(\phi;s)$ and
\[
\kappa_{m,n}=\frac{\pi^{mn}}{\prod_{j=0}^{m-1}j!(n-m+j)!}>0.
\]
Singular-value integration gives
\begin{equation}\label{eq:svd}
Z(s)=\kappa_{m,n}
\int_{\lambda_1>\cdots>\lambda_m>0}
e_r(\lambda)^s\Delta_m(\lambda)^2
\prod_{i=1}^m\lambda_i^{n+r-m}e^{-\sum_i\lambda_i}\,d\lambda,
\end{equation}
Write $\sigma_i=\sqrt{\lambda_i}$ and differentiate
$X=U_0(\operatorname{diag}(\sigma_1,\ldots,\sigma_m)\ 0)V_0^*$, with
$U_0,V_0$ unitary. At a diagonal point, each additional complex column
direction contributes $\sigma_i^2$ to the Jacobian, giving
$\sigma_i^{2(n-m)}$ in total. For $i<j$, the two independent
complex angular coordinates $a,b$ map, after conjugating one output
coordinate, to
\[
(a\sigma_j-\sigma_i b,\ -a\sigma_i+\sigma_j b).
\]
Their complex determinant is $\sigma_j^2-\sigma_i^2$, so their real
Jacobian is $(\lambda_i-\lambda_j)^2$. The diagonal phase and
radial coordinate contribute $\sigma_i\,d\sigma_i
=\tfrac12d\lambda_i$. Integration over the angular variables
therefore leaves a constant times
$\Delta_m(\lambda)^2\prod_i\lambda_i^{n-m}\,d\lambda$.
The weight $\det(XX^*)^r$ supplies $\prod_i\lambda_i^r$.

The constant follows from the unweighted Gaussian integral
$\pi^{mn}$. Expanding both Vandermonde determinants and integrating
their monomials on $(0,\infty)^m$ gives $m!$ times the moment
determinant. Division by the $m!$ ordered chambers leaves
\[
\det[(n-m+i+j)!]_{0\le i,j<m}
=\prod_{i=0}^{m-1}(n-m+i)!\prod_{j=0}^{m-1}j!.
\]
Indeed, after factoring $(n-m+i)!$ from row $i$, column $j$
is the monic degree-$j$ polynomial
$\prod_{a=1}^j(n-m+i+a)$ in $n-m+i$. Its determinant is the
Vandermonde determinant at consecutive integers. This proves the
stated value of $\kappa_{m,n}$.

Use ratio coordinates $t_1=\lambda_1$ and
$t_j=\lambda_j/\lambda_{j-1}$ for $2\le j\le m$.
Thus $t_1>0$, $0<t_j<1$, and $\lambda_j=t_1\cdots t_j$.
Put
\begin{align*}
a_k(s)&=(m-k)(n+r-k)+(r-k)_+s\quad(0\le k<m),\\
L&=1+t_2+t_2t_3+\cdots+t_2\cdots t_m,\\
D&=\prod_{1\le i<j\le m}(1-t_{i+1}\cdots t_j)^2,\qquad
U=\frac{e_r(\lambda)}{\lambda_1\cdots\lambda_r},
\end{align*}
where $(r-k)_+=\max(r-k,0)$. The last expression is a polynomial
in $t_2,\ldots,t_m$. For a term $\lambda_{i_1}\cdots\lambda_{i_r}$
with $i_1<\cdots<i_r$, one has $i_a\ge a$; hence its quotient by
$\lambda_1\cdots\lambda_r$ has nonnegative ratio exponents.
The leading term gives one. Thus on the closed cube,
$1\le U\le\binom mr$ and $1\le L\le m$.
The ratio-coordinate Jacobian is $\prod_{j=1}^m t_j^{m-j}$.
Combining it with the Vandermonde monomial and the density, the
exponent of $t_j$ before the factor $e_r^s$ is
\[
(m-j+1)(n+r-m)+(m-j)(m-j+1)+(m-j)
=(m-j+1)(n+r-j+1)-1.
\]
The factor $(\lambda_1\cdots\lambda_r)^s$ adds
$(r-j+1)_+s$.
Integrating $t_1$ in \eqref{eq:svd} gives
\begin{equation}\label{eq:ratiointegral}
Z(s)=\kappa_{m,n}\Gamma(a_0(s))
\int_{(0,1)^{m-1}}\prod_{j=2}^m t_j^{a_{j-1}(s)-1}
                         L^{-a_0(s)}DU^s\,dt.
\end{equation}

Set $p=r-1$ and $h=(m-r)n-r\ge n-r\ge1$.
For $0\le k<r$, substituting $\rho=(m-r+1)(n+1)$ and expanding gives
\[
-a_k(-\rho)=(r-k-1)\big((m-r)n-(r-k)\big).
\]
In terms of $p,h$, this is
\begin{equation}\label{eq:qk}
q_k:=-a_k(-\rho)=(p-k)(h+k)\quad(0\le k<r).
\end{equation}
In particular $q_p=0$, and $a_k(-\rho)>0$ for $k\ge r$.

We use Taylor subtraction to compute the highest Mellin pole.
For a function $H(t,s)$ smooth in $0\le t\le1$ and holomorphic in
$s$, choose a positive integer $M$ and write
\[
H(t,s)=\sum_{\ell=0}^{M-1}\frac{\partial_t^\ell H(0,s)}{\ell!}t^\ell
       +t^M R_M(t,s).
\]
For a holomorphic exponent $a(s)$ with $\Re a(s)+M>0$,
continuation of the Mellin integral is
\[
\int_0^1t^{a(s)-1}H(t,s)\,dt
=\sum_{\ell=0}^{M-1}
 \frac{\partial_t^\ell H(0,s)}{\ell!(a(s)+\ell)}
 +\int_0^1t^{a(s)+M-1}R_M(t,s)\,dt.
\]
The remainder is holomorphic. If $a(s_0)=-q$ with an integer
$0\le q<M$, only the term
$\ell=q$ has a pole. In \eqref{eq:ratiointegral}, the positive
lower bounds for $L,U$ give uniformly bounded derivatives of
$L^{-a_0(s)}DU^s$ on the closed cube for $s$ near $-\rho$.
Choose $M$ larger than all $q_k$ and apply this formula successively
to $t_2,\ldots,t_r$. Tail exponents remain positive, so dominated
integration applies to every remainder and coefficient.

Together with $\Gamma(a_0(s))$, there are exactly $r$ possible
simple pole factors at this point. A term has order $r$ only if it
uses the coefficient of degree $q_{j-1}$ in every $t_j$, $2\le j\le r$,
and the Gamma residue. Every other term has order at most $r-1$.
The coefficient of order $r$ is therefore obtained by evaluating
these analytic coefficients at $s=-\rho$.

Define the head factors
\[
D_p=\prod_{1\le i<j\le p}(1-t_{i+1}\cdots t_j)^2,
\qquad L_p=\sum_{i=1}^p\prod_{j=2}^i t_j,
\]
and the tail factors
\[
L_{\mathrm{tail}}=1+t_{r+1}+t_{r+1}t_{r+2}
                          +\cdots+t_{r+1}\cdots t_m,
\quad
D_{\mathrm{tail}}=\prod_{r\le i<j\le m}(1-t_{i+1}\cdots t_j)^2.
\]
Since $q_p=0$, taking the residue in $t_r$ evaluates the analytic
factor at $t_r=0$. Every factor of $D$ crossing from indices
$1,\ldots,p$ to $r,\ldots,m$ contains $t_r$ and becomes one.
In a term of $U$, the exponent of $t_r$ is the number of selected
indices at least $r$, minus one. It is zero precisely when all
first $p$ indices and one tail index are selected. Thus
\[
D|_{t_r=0}=D_pD_{\mathrm{tail}},\qquad
L|_{t_r=0}=L_p,\qquad U|_{t_r=0}=L_{\mathrm{tail}}.
\]
Since $-a_0(-\rho)=ph$, the head coefficient is
\begin{equation}\label{eq:headcoefficient}
K_{p,h}=\left[\prod_{k=1}^{p-1}t_{k+1}^{(p-k)(h+k)}\right]
 D_pL_p^{ph}.
\end{equation}
For $p=1$, both head factors and the coefficient in zero variables
equal one.

To calculate \eqref{eq:headcoefficient}, use variables
$z_1,\ldots,z_p$ with $z_i/z_1=t_2\cdots t_i$ for $i\ge2$, and set
$H=h+p-1$. The identity
\[
\Delta_p(z)^2=\left(\prod_{i=1}^p z_i^{2(p-i)}\right)D_p
\]
converts the selected ratio exponents into equal exponents $H$.
Indeed, in the homogeneous polynomial
$\Delta_p(z)^2(\sum_i z_i)^{ph}$ the exponent sum in positions
$k+1,\ldots,p$ is
\[
2\sum_{i=k+1}^p(p-i)+(p-k)(h+k)=(p-k)H.
\]
Successive differences force the exponents of $z_2,\ldots,z_p$
to equal $H$; total degree $pH$ forces the exponent of $z_1$ too.
Thus
\begin{equation}\label{eq:centralcoefficient}
K_{p,h}=[z_1^H\cdots z_p^H]\Delta_p(z)^2(\textstyle\sum_i z_i)^{ph}.
\end{equation}
Let $\mathfrak S_p$ be the permutation group of $\{0,\ldots,p-1\}$,
and let $\operatorname{sgn}(\sigma)$ denote a permutation's sign.
Set $1/a!=0$ for negative integers $a$.
Indexing determinant powers by $0,\ldots,p-1$ and expanding gives
\[
(ph)!\sum_{\sigma,\tau\in\mathfrak S_p}
 \operatorname{sgn}(\sigma)\operatorname{sgn}(\tau)
 \prod_{i=0}^{p-1}\frac1{(H-\sigma(i)-\tau(i))!}.
\]
Reindexing the product by $\sigma(i)$ leaves a determinant sum over
$\tau\sigma^{-1}$ for each of the $p!$ choices of $\sigma$.
Consequently,
\begin{align}
K_{p,h}
 &=(ph)!p!\det\left[\frac1{(H-i-j)!}\right]_{0\le i,j<p}\notag\\
 &=(-1)^{p(p-1)/2}p!(ph)!
                         \prod_{j=0}^{p-1}\frac{j!}{(h+j)!}\ne0.
\label{eq:factorialdet}
\end{align}
For the second equality, factor $1/(H-i)!$ from row $i$.
The remaining entry is the falling factorial
$(H-i)(H-i-1)\cdots(H-i-j+1)$, a monic polynomial of degree $j$
in $H-i$. Its determinant is the Vandermonde determinant of
$H,H-1,\ldots,H-p+1$, which has sign $(-1)^{p(p-1)/2}$ and
absolute value $\prod_{j=0}^{p-1}j!$.

The remaining tail integral is
\[
J=\int_{(0,1)^{m-r}}
\prod_{j=r+1}^m t_j^{a_{j-1}(-\rho)-1}
D_{\mathrm{tail}}L_{\mathrm{tail}}^{-\rho}
\,dt_{r+1}\cdots dt_m.
\]
Every exponent parameter in this integral is positive, the factor
$L_{\mathrm{tail}}$ is bounded away from zero, and the integrand is
positive in the open cube. Hence $0<J<\infty$.

The residue of $\Gamma(a_0(s))$ is
$(-1)^{ph}/((ph)!r)$, and the other critical slopes are
$r-1,r-2,\ldots,1$. Combining these with \eqref{eq:factorialdet},
\begin{equation}\label{eq:highestresidue}
\lim_{s\to-\rho}(s+\rho)^rZ(s)
=\frac{\kappa_{m,n}}{r}
(-1)^{ph+p(p-1)/2}
\left(\prod_{j=0}^{p-1}\frac{j!}{(h+j)!}\right)J\ne0.
\end{equation}
This proves the asserted pole order.
\end{proof}

\subsection{Determination of all integral multiplicities}\label{subsec:integermultiplicity}

\begin{proposition}\label{prop:integers}
For $1\le r<m\le n$ and every $\rho\in\mathcal R\cap\Z$,
\[
\mu(\rho)=|J(\rho)|,\qquad
\sum_{u\in\Z}\mu(u)=r.
\]
\end{proposition}
\begin{proof}
For $r\ge2$, Lemma~\ref{lem:polebound} and
Proposition~\ref{prop:integerpole} give
$\sum_{u\in\Z}\mu(u)\ge r$. For $r=1$, the unweighted integral is
\[
Z_{I_1}(e^{-\|X\|^2};s)
=\pi^{mn}\frac{\Gamma(s+mn)}{\Gamma(mn)},
\]
whose pole at $-mn$ gives the same bound.
For each $i$, the $r-i$ consecutive integers
$(m-i)(n-i)+u$, $0\le u<r-i$, contain exactly one multiple of $r-i$.
The upper bound \eqref{eq:upper} therefore contains exactly $r$
integral labels, counted with repetition. It follows that
\[
r\le\sum_{u\in\Z}\mu(u)
\le\sum_{\rho\in\mathcal R\cap\Z}|J(\rho)|=r.
\]
Each difference $|J(\rho)|-\mu(\rho)$ is nonnegative and their sum
is zero, so every difference is zero.
\end{proof}

\section{Multiplicities of the nonintegral roots}\label{sec:noninteger}
Throughout this section, assume $1\le r<m\le n$. We prove simplicity
by separating central supports in a filtration of the Fourier-transformed module. The needed facts about central
characters were recalled in Section~\ref{subsec:central}.

\subsection{The module generated by bounded words}\label{subsec:words}

Fix a nonintegral $\rho\in\mathcal R$. The polynomial in
\eqref{eq:upper} gives a finite primary decomposition of $Q$ for
$s$. Let $\delta_\rho$ be the projection of $\delta$ to the component
where $N_\rho=s+\rho$ is nilpotent, and now let $v$ be its Fourier
transform. All relations in \eqref{eq:fourierrelations} continue to
hold, and a fixed power of $N_\rho$ kills this whole component.
Put $z_0=m+n-\rho\notin\Z$, the value of $z$ at $s=-\rho$.
For $0\le k<r$, define $U(\g)$-submodules
\begin{equation}\label{eq:W}
W_k=U(\g)\C[s]\Span\{B^kv\},\qquad
W=\sum_{k=0}^{r-1}W_k.
\end{equation}
Using
$[A_{uv},x_{ij}]=\delta_{vi}x_{uj}$ and
$[A_{uv},\partial_{ij}]=-\delta_{ui}\partial_{vj}$ gives
\begin{align*}
\left[A_{uv},\sum_{i,j}x_{ij}\partial_{ib}\partial_{aj}\right]
={}&\sum_jx_{uj}\partial_{vb}\partial_{aj}
 -\sum_jx_{uj}\partial_{vb}\partial_{aj}\\
&-\delta_{ua}\sum_{i,j}x_{ij}\partial_{ib}\partial_{vj}.
\end{align*}
The first two sums cancel, and the term $z\partial_{ab}$ gives
\begin{equation}\label{eq:AB}
[A_{uv},B_{ab}]=-\delta_{ua}B_{vb}.
\end{equation}
Commuting a $B$-word past a product of row operators preserves its
length. It follows that $B_{ab}W_k\subseteq W_{k+1}$, with $W_r=0$,
and that every word of length $r$ kills $W$: each resulting term
contains at least $r$ commuting $B$ operators acting on $v$.
Thus $W$ is stable under $B_{ab}$, $\g$, and $s$.

Since $E$ commutes with $U(\g)$, the operator $E+rs+mn$ acts on
$W_k$ by the distinct scalar $-k$. Hence
\begin{equation}\label{eq:Wdirect}
W=\bigoplus_{k=0}^{r-1}W_k,\qquad E|_{W_k}=q-k.
\end{equation}
The space spanned by the derivative $(r-k)$-minors is stable under
commutators with $\g$: the identity
$[A_{uv},\partial_{ab}]=-\delta_{ua}\partial_{vb}$ replaces one
row in a derivative determinant. For such a minor $D$ and a product
$A_1\cdots A_t$ of row operators, repeatedly use
$DA_1=A_1D+[D,A_1]$. Every resulting rightmost constant-coefficient
operator is a linear combination of minors of the same size.
Equation~\eqref{eq:fourierrelations} kills each of them on $B^kv$.
Induction on $t$ therefore gives
\begin{equation}\label{eq:Cvanishing}
I_{r-k}(\partial)W_k=0,\qquad C_{r-k}W_k=0.
\end{equation}

Let $F_0=\bigcap_{a,b}\ker(B_{ab}:W\to W)$. The direct sum
\eqref{eq:Wdirect} is compatible with this kernel, since each
$B_{ab}$ sends $W_k$ to $W_{k+1}$. For $u\in F_0\cap W_k$,
\eqref{eq:commonkernel} and \eqref{eq:Cvanishing} imply
\begin{equation}\label{eq:F0polynomial}
\prod_{i=0}^{r-k-1}(q-k+iz-i^2)u=0.
\end{equation}
At $s=-\rho$, two different indices $i,j$ could give zero factors
only if $i(i-z_0)=j(j-z_0)$, which would force $z_0=i+j\in\Z$.
Each factor has nonzero $s$-slope $-(r-i)$. On the primary component
it is therefore either invertible or a nonzero scalar times
$N_\rho$. In the first case it has the form $c+dN_\rho$ with
$c\ne0$, whose inverse is the finite geometric series
$c^{-1}\sum_{a\ge0}(-dN_\rho/c)^a$; the sum terminates by nilpotence.
Equation~\eqref{eq:F0polynomial} proves
\begin{equation}\label{eq:NF0}
N_\rho F_0=0.
\end{equation}

When $F_0\cap W_k\ne0$, the unique vanishing index $i$ satisfies
$r\rho-mn-k=i(i-z_0)$. Equation~\eqref{eq:commonkernel} determines all its
central eigenvalues. They are those of the weight
\[
\lambda^{(i)}=(\underbrace{i-z_0,\ldots,i-z_0}_{i\text{ entries}},
                         \underbrace{0,\ldots,0}_{m-i\text{ entries}}),
\]
because
\begin{align}
c_p(\lambda^{(i)})
 &=\binom ip(i-z_0)(i-z_0+1)\cdots(i-z_0+p-1)\notag\\
 &=\frac1{p!}\prod_{a=0}^{p-1}\big(i(i-z_0)+az_0-a^2\big).
\label{eq:rectangle}
\end{align}
For $p>i$, both sides are zero. For $p\le i$, every surviving term
in \eqref{eq:capellieigenvalue} uses only the first $i$ coordinates,
which gives the first equality. The second uses
$i(i-z_0)+az_0-a^2=(i-a)(i+a-z_0)$.
Thus $F_0$ has finitely many exact central characters.

\subsection{Separation of the filtration layers}\label{subsec:separation}

Define $F_{-1}=0$ and, for $j\ge0$,
\[
F_j=\{w\in W:B_{\eta_1}\cdots B_{\eta_{j+1}}w=0
                     \text{ for all entry indices }\eta_1,\ldots,\eta_{j+1}\}.
\]
These are $U(\g)$-submodules, and $F_{r-1}=W$.
Take the formal word representation
$T=(\C^m)^*\otimes\C^n$, where the star denotes the linear dual,
$\g$ acts on the first factor, and the second is a multiplicity space.
Write $\Sym^jT$ for the $j$th symmetric tensor power of $T$.
By \eqref{eq:AB}, a basis
element of $T$ acts as the corresponding $B_{ab}$.
For $w\in F_j$, send a formal degree-$j$ word to the corresponding
$B$-word applied to $w$. Commutativity of the $B$ operators makes
this a linear map from $\Sym^jT$ to $F_0$, because any further
$B$ kills its value. The map from $w$ to this element of
$\operatorname{Hom}_\C(\Sym^jT,F_0)$ has kernel $F_{j-1}$.
The commutator relation \eqref{eq:AB} makes it equivariant: the
action of a row operator on a word is its commutator with that
word. Since $\Sym^jT$ is finite-dimensional, this gives an injection
\begin{equation}\label{eq:layerinjection}
F_j/F_{j-1}\lhook\joinrel\longrightarrow
                          (\Sym^jT)^*\otimes F_0.
\end{equation}
The weights of $(\Sym^jT)^*$ have nonnegative integral coordinates
$\nu=(\nu_1,\ldots,\nu_m)$ with $\sum_a\nu_a=j$.
Lemma~\ref{lem:tensorsupport} gives a finite-codimensional central
annihilator of each layer, with support contained in
\begin{equation}\label{eq:Sigma}
\Sigma_j=\{\chi_{\lambda^{(i)}+\nu}:0\le i<r,
                         \ \nu\in\Z_{\ge0}^m,\ |\nu|=j\}.
\end{equation}

\begin{lemma}\label{lem:separation}
If $j\ne\ell$, then $\Sigma_j\cap\Sigma_\ell=\varnothing$.
\end{lemma}
\begin{proof}
The shifted coordinates of $\lambda^{(i)}+\nu$ have exactly $i$
entries in the nonzero class $-z_0+\Z$ of $\C/\Z$, with the
remaining entries in $\Z$. Equality of central characters means
equality of these coordinate multisets, by \eqref{eq:HC}. Thus
$\chi_{\lambda^{(i)}+\nu}=\chi_{\lambda^{(i')}+\nu'}$ first forces
$i=i'$. Equality of the trace character then gives
$i(i-z_0)+j=i(i-z_0)+\ell$, hence $j=\ell$.
\end{proof}

\begin{proposition}\label{prop:nonintegers}
For $1\le r<m\le n$, every nonintegral root of $b_{I_r}(s)$ is simple.
\end{proposition}
\begin{proof}
For $w\in F_j$, every length-$j$ word $B^jw$ belongs to $F_0$.
Commutativity and \eqref{eq:NF0} give
$B^jN_\rho w=N_\rho B^jw=0$, hence
$N_\rho F_j\subseteq F_{j-1}$. We prove $N_\rho F_j=0$
by induction. The case $j=0$ is \eqref{eq:NF0}. If
$N_\rho F_{j-1}=0$, then $N_\rho$ induces a $Z$-linear map
\[
F_j/F_{j-1}\longrightarrow F_{j-1}.
\]
Set $J_j=\Ann_Z(F_j/F_{j-1})$. This ideal has finite codimension
and its zero set is the actual central support of the layer, hence
is contained in $\Sigma_j$. The product $J_0\cdots J_{j-1}$
annihilates $F_{j-1}$: its successive factors lower the filtration
until reaching $F_{-1}=0$. A finite product of finite-codimensional ideals in the
polynomial algebra $Z$ is again finite-codimensional, so this also
gives a finite central support for $F_{j-1}$. The zero set of a product is the union of the zero
sets of its factors. Lemma~\ref{lem:separation} makes the two zero
sets disjoint, so the Nullstellensatz gives
\[
J_j+J_0\cdots J_{j-1}=Z.
\]
Choose elements from the two ideals summing to one. Both kill the
image of the displayed map, so the image is zero. This completes
the induction and proves $N_\rho W=0$.

In particular, $(s+\rho)v=0$. After inverse Fourier transformation,
$(s+\rho)\delta_\rho=0$. The projected generator is cyclic over
$\gr_V^0D_{X,t}$ and $s$ is central there, so $s+\rho$ kills the
entire primary component of $Q$. Its exponent in the minimal
polynomial is therefore at most one.
\end{proof}

\section{Root existence and completion of the formula}\label{sec:existence}

\subsection{The known root-existence theorems}\label{subsec:knownroots}
We use the following result directly; its proof
is not repeated here. It is stated in the same unshifted ideal
normalization as in Section~\ref{subsec:graph}.

\begin{theorem}[Chen--Zuo, {\cite[Theorem~3.9]{CZ}}]
\label{prop:rootexistence}
For $1\le r<m\le n$, and integers $0\le i<r$ and $0\le u<r-i$,
\[
b_{I_r^{m,n}}\left(-\frac{(m-i)(n-i)+u}{r-i}\right)=0.
\]
Thus every element of $-\mathcal R$ occurs as a root.
\end{theorem}

For the boundary case we use the full result of
L\H orincz, Raicu, Walther, and Weyman:
\begin{theorem}[{\cite[Theorem~4.1]{LRWW}}]\label{thm:knownmaximal}
For $1\le m\le n$,
\[
b_{I_m^{m,n}}(s)=\prod_{j=n-m+1}^{n}(s+j).
\]
\end{theorem}
Their convention has the larger matrix dimension first; transposing
the generic matrix gives exactly the convention used here.

\subsection{Proof of the polynomial formula}\label{subsec:completion}
\begin{proof}[Proof of Theorem~\ref{thm:main}]
Suppose first that $r<m$. Proposition~\ref{prop:upper} excludes all
roots outside $-\mathcal R$, and Theorem~\ref{prop:rootexistence}
establishes every root in that set. Proposition~\ref{prop:integers}
gives the integral multiplicities, while
Proposition~\ref{prop:nonintegers} makes each nonintegral root simple.
This proves \eqref{eq:multiplicity}. The unique integral value in
the $i$th window is $\lceil(m-i)(n-i)/(r-i)\rceil$, which gives
\eqref{eq:mainproduct}. Theorem~\ref{thm:knownmaximal} proves
\eqref{eq:maximal}. For $r=1<m$, the sole label is $\rho=mn$;
for $r=m=1$, use Theorem~\ref{thm:knownmaximal}. Thus in every case
$b_{I_1^{m,n}}(s)=s+mn$.
\end{proof}

\begin{corollary}\label{cor:degree}
If $r<m$, then
\[
\sum_{\rho\in\Z}\mult_{s=-\rho}b_{I_r}(s)=r,\qquad
\deg b_{I_r}=r+|\mathcal R\setminus\Z|.
\]
For the determinantal variety $L_r=V(I_r)$, of codimension
$c=(m-r+1)(n-r+1)$, the embedding-independent convention of
\cite{BMS} is
\[
b_{L_r}(s)=b_{I_r}(s-c).
\]
Its roots are $c-\rho$, with the same multiplicities as in
Theorem~\ref{thm:main}.
\end{corollary}
\begin{proof}
The first two assertions follow by counting the factors in
\eqref{eq:mainproduct}. For $k>0$, choose an invertible $k\times k$ block and write
$X=\left(\begin{smallmatrix}A&B\\ C&D\end{smallmatrix}\right)$.
On the rank-$k$ locus, $A,B,C$ are free and $D=CA^{-1}B$.
The rank-$k$ stratum therefore has dimension
$k^2+k(n-k)+(m-k)k=k(m+n-k)$.
The rank-zero stratum is a point. These dimensions are increasing
for $0\le k<r\le m\le n$, so
$\dim L_r=(r-1)(m+n-r+1)$ and
$mn-\dim L_r=(m-r+1)(n-r+1)$.
The definition of the variety polynomial then gives the displayed
shift, which preserves every root multiplicity.
\end{proof}

\section{Multiplicity bounds and strong motivic monodromy}\label{sec:arithmetic}

The integral multiplicities admit a direct interval count, and this
count gives a sharp bound independent of the matrix dimensions.

\subsection{An interval formula for integral multiplicities}\label{subsec:interval}

\begin{corollary}\label{cor:interval}
Suppose $r<m\le n$, and set $C=(m-r)(n-r)$. For an integer $T$, let
\[
H_C(T;r)=\#\{d\in\Z:1\le d\le r,\ d^2-Td+C\le0\}.
\]
For an integer $\rho$, take $T=\rho-(m+n-2r)$. Then
\begin{equation}\label{eq:intervalcount}
\mult_{s=-\rho}b_{I_r}(s)=H_C(T;r)-H_C(T-1;r).
\end{equation}
If $T<0$ or $T^2<4C$, then $H_C(T;r)=0$. Otherwise,
\begin{equation}\label{eq:closedcount}
\begin{aligned}
H_C(T;r)=\max\Bigl\{0,\;&
\min\left(r,\left\lfloor\frac{T+\sqrt{T^2-4C}}{2}\right\rfloor\right)\\
&-\max\left(1,\left\lceil\frac{T-\sqrt{T^2-4C}}{2}\right\rceil\right)+1\Bigr\}.
\end{aligned}
\end{equation}
\end{corollary}
\begin{proof}
Substitute $d=r-i$ in the integral factors of
\eqref{eq:mainproduct}. Since
\[
\frac{(m-i)(n-i)}{r-i}=m+n-2r+d+\frac{C}{d},
\]
the relevant indices are exactly those with
$T-1<d+C/d\le T$. This proves \eqref{eq:intervalcount}.
For the closed formula, complete the square:
$d^2-Td+C\le0$ if and only if
$(2d-T)^2\le T^2-4C$. Counting the integers $d$ in $[1,r]$ that
satisfy $|2d-T|\le\sqrt{T^2-4C}$ gives \eqref{eq:closedcount}.
\end{proof}

\subsection{The sharp uniform multiplicity bound}\label{subsec:sharp}

\begin{proof}[Proof of Theorem~\ref{thm:sharpbound}]
Nonintegral roots and roots for maximal minors are simple. For the
other cases Corollary~\ref{cor:interval} reduces the problem to
bounding the size $N$ of
\[
S=\{d\in\Z:1\le d\le r,\ d+\lceil C/d\rceil=T\},\qquad C\ge1.
\]
If $S$ is nonempty, let $R=\max S$ and put
\[
u=R(T-R)-C,\qquad P=2R-T.
\]
Since $R\in S$, one has $0\le u\le R-1$.
For an integer $1\le d\le R$, membership in $S$ is
equivalent to $0\le d(T-d)-C\le d-1$. On putting $d=R-j$,
\[
d(T-d)-C=u+j(2R-T)-j^2=u+Pj-j^2.
\]
Its lower and upper bounds give exactly
\begin{equation}\label{eq:shell}
j(j-P)\le u,\qquad j(P+1-j)\le R-1-u,\qquad j\ge0.
\end{equation}
Adding the first two inequalities gives $j\le R-1$; thus these
inequalities already enforce $1\le d\le R$.

The first inequality has integer solutions $0\le j\le J$ for some
$J\ge0$. If $P<0$, then
$J^2\le J(J-P)\le u\le R-1$. Since $N\le J+1$ and
$\lfloor(J+1)^2/4\rfloor\le J^2$ (also for $J=0$), we obtain
$\lfloor N^2/4\rfloor\le R-1$.

Suppose $P\ge0$. Then $J=P+k$ for an integer $k\ge0$.
For a nonnegative integer $t$, put $F(t)=\lfloor(t+1)^2/4\rfloor$.
The largest integer value of $j(P+1-j)$ is $F(P)$.
If $R-1-u\ge F(P)$, the second inequality removes no point, so $N=J+1$.
Using $F(t+1)-F(t)=\lfloor(t+2)/2\rfloor\le J$ for $P\le t<J$,
\[
R-1\ge J(J-P)+F(P)\ge F(J)=\lfloor N^2/4\rfloor.
\]
If $R-1-u<F(P)$, the integer solutions of the second inequality are
$0\le j\le z$ or $j\ge P+1-z$, where $2z\le P$.
Their intersection with $[0,P+k]$ has size $N=2z+k+1$;
when $z=k=0$ the right interval contributes no point.
Since $u\ge k(P+k)$ and $R-1-u\ge z(P+1-z)$,
\begin{align*}
R-1&\ge k(P+k)+z(P+1-z)\\
   &\ge(z+k)^2+z
    \ge\left\lfloor\frac{(2z+k+1)^2}{4}\right\rfloor.
\end{align*}
The difference in the second step is $(k+z)(P-2z)\ge0$.
The difference in the last step is
$zk+k^2-\lfloor(k+1)^2/4\rfloor\ge0$, including $k=0$.
In every case $\lfloor N^2/4\rfloor\le R-1$, so
$N^2<4R\le4r$. This proves the upper bound in \eqref{eq:sharpbound}.

To attain it, set $L=\lfloor\sqrt{4r-1}\rfloor$ and choose
\[
C=r(r-L+1),\qquad T=2r-L+1.
\]
The inequality $L<2\sqrt r\le r+1$ gives $L\le r$, hence $C\ge1$.
Moreover,
\[
d^2-Td+C=(d-r)(d-(r-L+1)),\qquad
(T-1)^2-4C=L^2-4r<0.
\]
Thus $H_C(T;r)=L$ and $H_C(T-1;r)=0$.
The dimensions $m=r+1$, $n=r+C$ realize this value of $C$ and
satisfy $r<m\le n$. Equation~\eqref{eq:intervalcount} gives a root
of multiplicity $L$.
\end{proof}

\subsection{Strong motivic monodromy with multiplicities}\label{subsec:strong}
We now compare the denominator supplied by the standard determinantal
resolution with the polynomial in Theorem~\ref{thm:main}.

\begin{proposition}[Numerical data of the determinantal resolution]
\label{prop:detresolution}
There is a log resolution $\pi:Y\to\Mat_{m,n}$ of $I_r$ whose
relevant divisors are $E_0,\ldots,E_{r-1}$, and
\begin{equation}\label{eq:detresolution}
I_r\mathcal O_Y=\mathcal O_Y\left(-\sum_{i=0}^{r-1}(r-i)E_i\right),
\qquad
K_{Y/\Mat_{m,n}}=\sum_{i=0}^{r-1}\big((m-i)(n-i)-1\big)E_i.
\end{equation}
\end{proposition}
This is the standard resolution by successive blow-ups of the strict
transforms of the loci of rank at most $0,1,\ldots,r-1$; see
\cite{Vainsencher} and \cite[Proposition~3.3]{CZ25} for the numerical
data. Our index $i$ is one less than the index in the latter reference.
When $r=m=n$, the final rank locus is a divisor, so its blow-up is
an isomorphism; $E_{m-1}$ denotes its strict transform and is included
in the total transform. Its data are $(N_{m-1},\nu_{m-1})=(1,1)$.
In particular, it must not be omitted merely because it is not
exceptional.

\begin{lemma}\label{lem:coincidentslopes}
Let $0\le i<j<r$. If
\[
\frac{(m-i)(n-i)}{r-i}
=\frac{(m-j)(n-j)}{r-j}=\rho,
\]
then
\begin{equation}\label{eq:sloperelation}
\rho=m+n-i-j\in\Z,
\qquad (m-r)(n-r)=(r-i)(r-j).
\end{equation}
Thus every nonintegral candidate slope occurs at only one index.
\end{lemma}
\begin{proof}
The two equalities say that $i$ and $j$ are roots of
\[
x^2+(\rho-m-n)x+mn-r\rho=0.
\]
Their sum gives $i+j=m+n-\rho$. Substituting this into the
product relation $ij=mn-r\rho$ gives the second identity in
\eqref{eq:sloperelation}.
\end{proof}

\begin{proof}[Proof of Theorem~\ref{thm:strong}]
Insert \eqref{eq:detresolution} into the resolution formula
\eqref{eq:motivicresolution}. In each summand an index occurs at
most once. Hence $D_r(T)$ in \eqref{eq:motivicdenominator} is a
common denominator and
\[
D_r(T)Z_{I_r}^{\mathrm{mot}}(T)\in\mathcal M_{\C}[T].
\]
In fact, the numerator is obtained by multiplying the $J$th
summand by the missing factors indexed by $i\notin J$.
Every denominator has slope $\nu_i/N_i$, and a fixed slope $\rho$
occurs at most $p_r(\rho)$ times. Definition~\ref{def:motivicorder}
therefore gives the first inequality in \eqref{eq:stronginequality}.

For the second inequality, first assume $r<m$. If $\rho\in\Z$,
each index counted by $p_r(\rho)$ supplies a distinct label $(i,0)$
in $J(\rho)$. Theorem~\ref{thm:main} gives
\[
p_r(\rho)\le |J(\rho)|=\mult_{s=-\rho}b_{I_r}(s).
\]
If $\rho\notin\Z$ and $p_r(\rho)>0$, then $\rho\in\mathcal R$;
Lemma~\ref{lem:coincidentslopes} and Theorem~\ref{thm:main} give
$p_r(\rho)=1=\mult_{s=-\rho}b_{I_r}(s)$. If $p_r(\rho)=0$,
the inequality is automatic.
For $r=m$, the slopes are the distinct integers $n-i$
($0\le i<m$), and each corresponding Bernstein--Sato root is
simple by \eqref{eq:maximal}. This proves the second inequality
also in the maximal-minor case.

Finally, set
\[
Q_r(s)=\prod_{i=0}^{r-1}(s+\nu_i/N_i).
\]
The preceding multiplicity comparison proves $Q_r(s)\mid b_{I_r}(s)$.
The denominator of every summand of
\eqref{eq:topologicaldefinition} divides $Q_r(s)$ up to a nonzero
rational scalar. Therefore $Q_r(s)Z_{I_r}^{\mathrm{top}}(s)$ is a
polynomial, and so is $b_{I_r}(s)Z_{I_r}^{\mathrm{top}}(s)$.
\end{proof}

\begin{remark}
The proof bounds the orders of all candidate slopes before any
cancellation. It does not require a separate assertion that every
candidate is an actual motivic pole. The motivic statement is the
simultaneous denominator bound, and the topological statement is the
ordinary rational-function divisibility
\eqref{eq:topologicalstrong}.
\end{remark}

\phantomsection

\end{document}